# Vaccination rules for a true-mass action SEIR epidemic model based on an observer synthesis. Preliminary results

M. De la Sen*, A. Ibeas** and S. Alonso-Quesada*

*Institute of Research and Development of Processes. Faculty of Science and Technology
University of the Basque Country. PO. Box 644- Bilbao. Spain
**Department of Telecommunications and Systems Engineering
Autonomous University of Barcelona, 08193- Bellaterra, Barcelona, Spain

**Abstract**. This paper presents a simple continuous-time linear vaccination-based control strategy for a SEIR (susceptible plus infected plus infectious plus removed populations) propagation disease model. The model takes into account the total population amounts as a refrain for the illness transmission since its increase makes more difficult contacts among susceptible and infected. The control objective is the asymptotically tracking of the removed-by-immunity population to the total population while achieving simultaneously the remaining population (i.e. susceptible plus infected plus infectious) to asymptotically converge to zero. A state observer is used to estimate the true various partial populations of susceptible, infected, infectious and immune which are assumed to be unknown. The model parameters are also assumed to be, in general, unknown. In this case, the parameters are replaced by available estimates to implement the vaccination action.

## 1. Introduction

Important control problems nowadays related to Life Sciences are the control of ecological models like, for instance, those of population evolution (Beverton-Holt model, Hassell model, Ricker model etc.) via the online adjustment of the species environment carrying capacity, that of the population growth or that of the regulated harvesting quota as well as the disease propagation via vaccination control. In a set of papers, several variants and generalizations of the Beverton-Holt model (standard time–invariant, time-varying parameterized, generalized model or modified generalized model) have been investigated at the levels of stability, cycle- oscillatory behavior, permanence and control through the manipulation of the carrying capacity (see, for instance, [1-5]). The design of related control actions has been proved to be important in those papers at the levels, for instance, of aquaculture exploitation or plague fighting. On the other hand, the literature about epidemic mathematical models is exhaustive in many books and papers. A non-exhaustive list of references is given in this manuscript, cf. [6-14] (see also the references listed therein). The various classes of epidemic models have two major variants, namely, the so-called "pseudo-mass action models", where the total population is not taken into account as a relevant disease contagious factor and the so-called "true-mass action models", where the total population is more realistically considered as an inverse factor of the disease transmission rates. There are many variants of the above models, for instance, including vaccination of different kinds as follows: constant [8-9], impulsive [9, 12], discrete – time etc., incorporating point or distributed delays [12-13], oscillatory behaviours [14] etc. . On the other hand, variants of such models become considerably simpler for the illness transmission among plants [6-7]. Recent work on the influence of latent and infective periods and studiers of persistence of infection have been performed in [20-21]. Other mathematical models involving coupled partial populations of the whole ones require very close analysis of equilibrium points, stability and positivity to that required in epidemic models (see, for instance, [22-23]).

In this paper, a continuous-time vaccination observer-based control strategy is given for a SEIR epidemic model which takes directly the estimated susceptible, infected, infectious and immune populations to design the vaccination strategy. It is not required either the knowledge through time of the true partial populations of susceptible, infected, infectious and immune, since this task is performed by the observer, nor the knowledge of the true parameters. The incorporation of the observer to the vaccination control rule is the main contribution of this paper. On the other hand, it is assumed that the total population is known and equal to the total estimated population and also that it remains constant through time, so that the illness transmission is not critical, and the SEIR – model is of the above mentioned true-mass action type. Some important issues of positivity, stability and model reference tracking of a suitable population evolution of the combined SEIR-model and its observer are discussed.

**1.1 Notation.** $R_+^n$ is the first open n-real orthant and $R_{0+}^n$ is the first closed n- real orthant.

$x \in R_{0+}^n$ is a positive real n-vector in the usual sense that all its components are nonnegative. This can be also denoted by $x > 0$ if $x \neq 0$. If all the components are strictly positive, this is denoted by $x >> 0$ or $x \in R_+^n$. In the same way, $A \in R_{0+}^{n \times n}$ is a positive real n-matrix in the usual sense that all its entries are nonnegative and we can use the notations $A > 0$ and $A >> 0$ (if $A \in R_+^{n \times n}$). A square real matrix A is a Metzler matrix if and only if all its off-diagonal entries are nonnegative and then its associate exponential matrix function is positive; i.e. none of its entries takes a negative value at any time.

$C^{(q)}(Do;Im)$ is the set of real functions of class q of domain Do and image Im. $PC^{(q)}(Do;Im)$ is the set of real functions of class (q-1) of domain Do and image Im whose q-th derivative exits but it is not necessarily everywhere continuous on its definition domain.

**2. SEIR epidemic model**

Let S (t) be the "susceptible" population of infection at time t, E (t) the " infected" ( i.e. those which incubate the illness but do not still have any symptoms) at time t, I (t ) is the " infectious" ( or "infective") population at time t, and R (t) is the " removed –by- immunity " ( or " immune") population at time t. Consider the true-mass action SEIR-type epidemic model:

$$\dot{S}(t) = -\mu S(t) + \omega R(t) - \beta \frac{S(t)I(t)}{N} + \mu N(1 - V(t)) \quad (1)$$

$$\dot{E}(t) = \beta \frac{S(t)I(t)}{N} - (\mu + \sigma)E(t) \quad (2)$$

$$\dot{I}(t) = -(\mu + \gamma)I(t) + \sigma E(t) \quad (3)$$

$$\dot{R}(t) = -(\mu + \omega)R(t) + \gamma I(t) + \mu N V(t) \quad (4)$$

subject to initial conditions $S_0 = S(0) \geq 0$, $E_0 = E(0) \geq 0$, $I_0 = I(0) \geq 0$ and $R_0 = R(0) \geq 0$ under the constraint $N = N(0) = S(t) + E(t) + I(t) + R(t) = S(0) + E(0) + I(0) + R(0)$ ; $\forall t \in R_{0+}$ and the vaccination function $V: R_{0+} \to R_{0+}$. In the above SEIR – model, N is the total population, $\mu$ is the rate of deaths from causes unrelated to the infection, ω is the rate of losing immunity, β is the transmission constant ( with the total number of infections per unity of time at time t being $\beta \frac{S(t)I(t)}{N}$), $\sigma^{-1}$ and $\gamma^{-1}$ are finite and , respectively, the average durations of the latent and infective periods. All the above parameters are assumed to be nonnegative. Note that $S, E, I, R \in C^{(1)}(\mathbf{R}_{0+}; \mathbf{R})$ if $V \in C^{(0)}(\mathbf{R}_{0+}; \mathbf{R})$. However, if $V \in PC^{(0)}(\mathbf{R}_{0+}; \mathbf{R})$ then $S, E, I, R \in C^{(0)}(\mathbf{R}_{0+}; \mathbf{R})$, their time-derivatives exist but they are not everywhere continuous on $\mathbf{R}_{0+}$, in general. This model has been studied in [19] from the point of view of equilibrium points in the free- vaccination case and control. Two vaccination auxiliary controls (grouping the last three terms of the right-hand-side of (1)) being respectively proportional to the susceptible or to the whole population so that the whole population is asymptotically immune have been proposed. Note that controllability is a very suitable property of practical control problems since each state variable is allowed to track a prescribed value, in a finite time (see, for instance, [24]). However, we have to point out that epidemic models are not , in general, controllable since all the partial populations cannot be separately controlled through vaccination since there is always an inherent interchange among the various populations, [19]. The SEIR model (1)-(4) is uncontrollable according to the above principle so that there are potential final values of some of the partial populations which are not within the set of reachable states, [25]. The epidemic models proposed in previous works assume that the parameters of the SEIR model are known. However, the research in this paper does not require such a perfect parametrical knowledge since an observer with parametrical estimation is incorporated.

**3. Observer-based vaccination strategy**

It turns out that while the assumption of the knowledge of the total population N is not quite restrictive in practice , the knowledge of the individual various partial populations of susceptible, infected, infectious and immune may be considered severely restrictive. It is seen that if the partial initial populations are unknown then their evolution through time cannot be computed in a closed form from the differential system (1)-(4). A practical solution to circumvent the problem might be to estimate them based on percentages of the total population through time from experimental knowledge of the disease propagation. Another solution may be to estimate them online by using an on-line observer. The current manuscript focuses this solution by using a SEIR – estimation algorithm (observer) of the SEIR model (1)-(4) which estimates through time the individual populations being involved. The vaccination strategy is obtained as a control strategy from the data supplied by the observer through time. Such a strategy does not require the knowledge of the partial populations to organize and perform the vaccination strategy. The estimates of the various individual populations are denoted by the same notations as the real populations with hat superscripts, namely, $\hat{S}(t), \hat{E}(t), \hat{I}(t), \hat{R}(t)$, Thus, consider the SEIR-type observer for the SEIR-epidemic model (1)-(4) as follows:

$$\dot{\hat{S}}(t) = -\hat{\mu}\hat{S}(t) + \hat{\omega}\hat{R}(t) - \hat{\beta}\frac{\hat{S}(t)\hat{I}(t)}{N} + \hat{\mu}N(1 - V(t)) \qquad (5)$$

$$\dot{\hat{E}}(t) = \hat{\beta}\frac{\hat{S}(t)\hat{I}(t)}{N} - (\hat{\mu} + \hat{\sigma})\hat{E}(t) \qquad (6)$$

$$\dot{\hat{I}}(t) = -(\hat{\mu} + \hat{\gamma})\hat{I}(t) + \hat{\sigma}\hat{E}(t) \qquad (7)$$

$$\dot{\hat{R}}(t) = -(\hat{\mu} + \hat{\omega})\hat{R}(t) + \hat{\gamma}\hat{I}(t) + \hat{\mu}NV(t) \qquad (8)$$

subject to initial conditions $\hat{S}_0 = \hat{S}(0) \geq 0$, $\hat{E}_0 = \hat{E}(0) \geq 0$, $\hat{I}_0 = \hat{I}(0) \geq 0$ and $\hat{R}_0 = \hat{R}(0) \geq 0$ under the constant population constraint

$$N = N(0) = \hat{S}(t) + \hat{E}(t) + \hat{I}(t) + \hat{R}(t) = \hat{S}(0) + \hat{E}(0) + \hat{I}(0) + \hat{R}(0)$$

equalizing its total estimated value for all time; $\forall t \in \mathbf{R}_{0+}$ and the vaccination law $V: \mathbf{R}_{0+} \to \mathbf{R}_{0+}$ generated by

$$V(t) = \frac{1}{\hat{\mu}N}\left(k_1\hat{S}(t) + k_2\hat{E}(t) + k_3\hat{I}(t) + k_4\hat{R}(t) + k_5\hat{S}(t)\hat{I}(t) + gN\right) \qquad (9)$$

where $k_i$; i=1,2, …, 5 are constant real gains to be determined in order to achieve the control objective. In the above estimated SEIR– model, $\hat{\mu}$ is the estimate of $\mu$, i.e. the estimated rate of deaths from causes unrelated to the infection and also the estimated rate of births, $\hat{\omega}$ is the estimate of $\omega$, i.e. the estimated rate of losing immunity, $\hat{\beta}$ is the estimate of $\beta$, i.e. the estimated transmission constant (with the total number of infections per unity of time at time t being $\hat{\beta}\frac{\hat{S}(t)\hat{I}(t)}{N}$), $\hat{\sigma}^{-1}$ and $\hat{\gamma}^{-1}$ are finite and, respectively, the estimates of $\sigma^{-1}$ and $\gamma^{-1}$, i.e. the estimated average durations of the latent and infective periods, respectively. The estimations of the above parameters can be done, through the use of available "a priori" knowledge, to be identical to the true values if those ones are known or estimated on-line from data measurements. Through this manuscript, we assume that those estimated parameters are fixed but not necessarily identical to the true parameters and all of them are nonnegative. The substitution of (9) in (5)-(8) yields the following combined observer- controller for the SEIR model:

$$\dot{\hat{x}}(t) = \hat{A}(t)\hat{x}(t) + \hat{b} \qquad (10)$$

where:

$$\hat{x}(t) := (\hat{S}(t), \hat{E}(t), \hat{I}(t), \hat{R}(t))^T \quad ; \quad \hat{b} := ((\hat{\mu} - g)N, 0, 0, gN)^T \qquad (11)$$

$$\hat{A}(t) := \begin{bmatrix} -(\hat{\mu} + k_1 + (\hat{\beta}_1 + k_5)\hat{I}(t)) & -k_2 & -k_3 & \hat{\omega} - k_4 \\ \hat{\beta}_1\hat{I}(t) & -(\hat{\mu} + \hat{\sigma}) & 0 & 0 \\ 0 & \hat{\sigma} & -(\hat{\mu} + \hat{\gamma}) & 0 \\ k_1 + k_5\hat{I}(t) & k_2 & \hat{\gamma} + k_3 & -(\hat{\mu} + \hat{\omega} - k_4) \end{bmatrix} ; \hat{\beta}_1 := \hat{\beta}/N \qquad (12)$$

The substitution of (9) into (1)-(4) yields the following SEIR observer-based vaccination controlled SEIR- model:

$$\dot{x}(t) = A(t)x(t) + B(t)\hat{x}(t) + b \qquad (13)$$

where:

$$x(t):=(S(t),E(t),I(t),R(t))^T \quad ; \quad b:=\left(\left(1-\frac{g}{\hat{\mu}}\right)\mu N,0,0,\frac{g\mu N}{\hat{\mu}}\right)^T \tag{14}$$

$$A(t):=\begin{bmatrix} -(\mu+\beta_1 I(t)) & 0 & 0 & \omega \\ \beta_1 I(t) & -(\mu+\sigma) & 0 & 0 \\ 0 & \sigma & -(\mu+\gamma) & 0 \\ 0 & 0 & \gamma & -(\mu+\omega) \end{bmatrix} \quad ; \quad \beta_1:=\beta/N \tag{15}$$

$$B(t):=\left(\frac{\mu}{\hat{\mu}}\right)\begin{bmatrix} -(k_1+k_5\hat{I}(t)) & -k_2 & -k_3 & -k_4 \\ 0 & 0 & 0 & 0 \\ 0 & 0 & 0 & 0 \\ k_1+k_5\hat{I}(t) & k_2 & k_3 & k_4 \end{bmatrix} \tag{16}$$

The systems (10)-(12) and (13)-(16) may be compacted as an extended system as follows:

$$\dot{\bar{x}}(t):=\bar{A}(t)\bar{x}(t)+\bar{b} \tag{17}$$

with

$$\bar{A}(t):=\begin{bmatrix} \hat{A}(t) & 0 \\ A(t)-\hat{A}(t)+B(t) & A(t) \end{bmatrix} \quad ; \quad \bar{b}:=\left(\hat{b}^T,\tilde{b}^T\right)^T \tag{18}$$

$\bar{x}(t):=\left(\hat{x}^T(t),\tilde{x}^T(t)\right)^T$ with $\tilde{x}(t)=x(t)-\hat{x}(t)$ being the observation error and $\tilde{b}$ is a parametrical error defined by

$$\tilde{b}:=b-\hat{b}=\left(\frac{\mu}{\hat{\mu}}-1\right)N(\hat{\mu}-g,0,0,g)^T \tag{19}$$

It is direct to see that $\|\tilde{b}\|\leq\varepsilon$ for any given real $\varepsilon>0$, with $\|\tilde{b}\|:=\left|\frac{\mu}{\hat{\mu}}-1\right|N\sqrt{(\hat{\mu}-g)^2+g^2}$ being the Euclidean norm of $\tilde{b}$, if $|\mu-\hat{\mu}|\leq\frac{\hat{\mu}\varepsilon}{N\sqrt{(\hat{\mu}-g)^2+g^2}}$. Decompose

$$A(t):=A_0+\Delta A(t) \quad ; \quad \hat{A}(t):=\hat{A}_0+\Delta\hat{A}(t) \quad ; \quad A(t)-\hat{A}(t)+B(t)=B_0+\Delta B(t) \tag{20}$$

where $A_0$, $\hat{A}_0$ and $B_0$ are constant matrices being the non-unique decompositions (20) as follows:

$$A_0:=\begin{bmatrix} -(\mu+\beta_1 I_r) & 0 & 0 & \omega \\ 0 & -(\mu+\sigma) & 0 & 0 \\ 0 & \sigma & -(\mu+\gamma) & 0 \\ 0 & 0 & \gamma & -(\mu+\omega) \end{bmatrix} \quad ; \quad \Delta A(t):=\begin{bmatrix} \beta_1(I_r-I(t)) & 0 & 0 & 0 \\ \beta_1 I(t) & 0 & 0 & 0 \\ 0 & 0 & 0 & 0 \\ 0 & 0 & 0 & 0 \end{bmatrix}$$

$$\hat{A}_0:=\begin{bmatrix} -(\hat{\mu}+k_1+(\hat{\beta}_1+k_5)\hat{I}_r) & -k_2 & -k_3 & \hat{\omega}-k_4 \\ 0 & -(\hat{\mu}+\hat{\sigma}) & 0 & 0 \\ 0 & \hat{\sigma} & -(\hat{\mu}+\hat{\gamma}) & 0 \\ k_1+k_5\hat{I}_r & k_2 & \hat{\gamma}+k_3 & -(\hat{\mu}+\hat{\omega}-k_4) \end{bmatrix} \quad ; \quad \Delta\hat{A}(t):=\begin{bmatrix} (\hat{\beta}_1+k_5)(\hat{I}_r-\hat{I}(t)) & 0 & 0 & 0 \\ \hat{\beta}_1\hat{I}(t) & 0 & 0 & 0 \\ 0 & 0 & 0 & 0 \\ k_5(\hat{I}(t)-\hat{I}_r) & 0 & 0 & 0 \end{bmatrix}$$

$$\Delta B(t)=A(t)-\hat{A}(t)+B(t)-B_0$$

$$
= \begin{bmatrix} \hat{\mu} - \mu + \hat{\beta}_1 \hat{I}(t) - \beta_1 I(t) + \left(1 - \frac{\mu}{\hat{\mu}}\right)(k_1 + k_5 \hat{I}(t)) - b_{011} & \left(1 - \frac{\mu}{\hat{\mu}}\right)k_2 & \left(1 - \frac{\mu}{\hat{\mu}}\right)k_3 & \omega - \hat{\omega} + \left(1 - \frac{\mu}{\hat{\mu}}\right)k_4 \\ \beta_1 I(t) - \hat{\beta}_1 \hat{I}(t) - b_{021} & \hat{\mu} - \mu + \hat{\sigma} - \sigma & 0 & 0 \\ 0 & \sigma - \hat{\sigma} & \hat{\mu} + \hat{\gamma} - \mu - \gamma & 0 \\ -\left(1 - \frac{\mu}{\hat{\mu}}\right)(k_1 + k_5 \hat{I}(t)) & -\left(1 - \frac{\mu}{\hat{\mu}}\right)k_2 & \gamma - \hat{\gamma} - \left(1 - \frac{\mu}{\hat{\mu}}\right)k_3 & \hat{\mu} + \hat{\omega} - \mu - \omega - \left(1 - \frac{\mu}{\hat{\mu}}\right)k_4 \end{bmatrix}
$$

(21)

for any given prefixed constant values $I_r \geq 0$, $\hat{I}_r \geq 0$ and $B_0 := \begin{bmatrix} b_{011} & 0 & 0 & 0 \\ b_{021} & 0 & 0 & 0 \\ 0 & 0 & 0 & 0 \\ 0 & 0 & 0 & 0 \end{bmatrix}$ so that $\overline{A}(t)$ is decomposed as follows:

$$\overline{A}(t) = \overline{A}_0 + \widetilde{\overline{A}}_0(t); \quad \overline{A}_0 := \begin{bmatrix} \hat{A}_0 & 0 \\ B_0 & A_0 \end{bmatrix}; \quad \widetilde{\overline{A}}_0(t) := \begin{bmatrix} \Delta \hat{A}(t) & 0 \\ \Delta B(t) & \Delta A(t) \end{bmatrix} \qquad (22)$$

If $A_0$ and $\hat{A}_0$ are stability (or Hurwitz) matrices then the block triangular matrix $\overline{A}_0$ is a stability matrix with stability abscissa $(-\rho)$ which is subject to $max\left(Re\,\lambda_i(A_0), Re\,\lambda_i(\hat{A}_0)\right) \leq -\rho < 0$ where the first inequality is non strict if there is some multiple eigenvalue of $\overline{A}_0$.

## 4. About the stability of the SEIR extended epidemic model (17)-(22)

The SEIR epidemic model (1)-(4) and its observer (10)-(12) are both always globally stable if $N < \infty$ if the initial real and observed populations are bounded and equalize N provided that both systems are positive as, on the other hand, the real problem at hand requires. This would imply that all the particular populations (susceptible, infected, infectious and immune) are nonnegative and less than N so that they are uniformly bounded for all time. Note that if any of the populations of the SEIR model is negative then boundedness of all the components of the associated system is not guaranteed. The boundedness of all the partial populations upper-bounded by a bounded totals population N guarantees that the extended SEIR system (17)-(22), constructed with the estimated and observer error states, is stable with $\sum_{i=1}^{4} \hat{x}_i(t) = N$ and Euclidean norm $\|\widetilde{x}(t)\| \leq \|x(t)\| + \|\hat{x}(t)\| \leq 2\sqrt{2}\,N$. This fact does not guarantee by itself that $\overline{A}(t)$ is a stability matrix function, that the observation error converges asymptotically to zero or that, at least, it converges to a small residual set around zero. This feature merits further discussion. Note that simple inspection shows that $A_0$ is a stability matrix if $d(s) = (s + \mu + \beta_1 I_r)(s + \mu + \sigma)(s + \mu + \gamma)(s + \mu + \omega)$ is Hurwitz. Also, $\hat{A}_0$ is a stability matrix if $det(sI - \hat{A}_0) = \hat{d}(s) + (k_1 + k_5 \hat{I}_r)\hat{n}(s)$ has all its zeros in $Re\,s < 0$ where

$$\hat{n}(s) = (k_4 - \hat{\omega})(s + \hat{\mu} + \hat{\sigma})(s + \hat{\mu} + \hat{\gamma}) \qquad (23)$$

$$\hat{d}(s) = (s + \hat{\mu} + k_1 + (\hat{\beta}_1 + k_5)\hat{I}_r)(s + \hat{\mu} + \hat{\sigma})(s + \hat{\mu} + \hat{\gamma})(s + \hat{\mu} + \hat{\omega} - k_4) \qquad (24)$$

Assume that $\hat{d}(s)$ is a Hurwitz polynomial, that is,

$$k_4 < \hat{\mu} + \hat{\omega} \quad ; \quad \hat{\mu} + k_1 + (\hat{\beta}_1 + k_5)\hat{I}_r > 0 \quad ; \quad \hat{\mu} + \hat{\sigma} > 0 \quad ; \quad \hat{\mu} + \hat{\gamma} > 0 \tag{25}$$

and define $\hat{h}(s) := (k_1 + k_5 \hat{I}_r)\hat{n}(s)/\hat{d}(s)$. Note that

$$det(s\,\mathrm{I} - \hat{A}_0) = \hat{d}(s) + (k_1 + k_5 \hat{I}_r)\hat{n}(s) = \hat{d}(s)(1 + \hat{h}(s)) \Leftrightarrow 1 + \hat{h}(s) = 0$$

has all its solutions in $\mathrm{Re}\,s < 0$ for all $\hat{n}(s)$ of the form (23) from the small gain theorem if and only if $\|\hat{h}\|_\infty := \sup_{\omega \in R_{0+}} \left|\frac{\hat{n}(i\omega)}{\hat{d}(i\omega)}\right| < 1$ ($i = \sqrt{-1}$ being the imaginary complex unit) since $\hat{d}(s)$ is a Hurwitz polynomial and, where $\|\hat{h}\|_\infty$ is the $H_\infty$-norm of the transfer function $\hat{h}(s)$. Since $\overline{A}_0$ is block-triangular and constant then the following result is direct.

**Assertion 1**. $\overline{A}_0$ is a stability matrix if and only if $\mu + \beta_1 I_r > 0$, $\mu + \sigma > 0$, $\mu + \gamma > 0$, $\mu + \omega > 0$ and $\hat{h} \in RH_\infty$ (i.e. $k_4 < \hat{\mu} + \hat{\omega}$ ; $\hat{\mu} + k_1 + (\hat{\beta}_1 + k_5)\hat{I}_r > 0$ ; $\hat{\mu} + \hat{\sigma} > 0$ ; $\hat{\mu} + \hat{\gamma} > 0$) with $\|\hat{h}\|_\infty < 1$. □

From Assertion 1 and Gronwall´s Lemma [17] the following follows:

**Assertion 2**. The matrix function $\overline{A}(t)$ is stable if $\overline{A}_0$ is a stability matrix, that is, if (25) holds and, furthermore, $\rho > \sup_{t \in R_{0+}} \|\tilde{\overline{A}}_0(t)\|$ where $(-\rho) < 0$ is the stability abscissa of the (stability) matrix $\overline{A}_0$.

□

Another alternative sufficiency-type stability condition of $\overline{A}(t)$, which replaces the constraint $\rho > \sup_{t \in R_{0+}} \|\tilde{\overline{A}}_0(t)\|$ in Assertion 2 (see [18] for close related problems) is as follows:

**Assertion 3**. The matrix function $\overline{A}(t)$ is stable, so that the unforced extended system (17)-(22) is then globally exponentially stable, if $\overline{A}_0$ is a stability matrix (see Assertion 1) and, furthermore, there exist norm-dependent real constants $\alpha_0 \in R_{0+}$ (being sufficiently small) and $\alpha_1 \in R_{0+}$ such that $\int_t^{t+T} \|\dot{\tilde{\overline{A}}}_0(\tau)\| d\tau \leq \alpha_0 T + \alpha_1$ for any given fixed $T \in R_+$ ; $\forall t \in R_{0+}$. □

Note that the Euclidean norm of $\overline{b}$ may be directly calculated from those of $\hat{b}$ and $\tilde{b}$ using (11) and (19) leading to:

$$\|\overline{b}\| \leq \frac{(\hat{\mu} + |\mu - \hat{\mu}|)N}{\hat{\mu}} \sqrt{(\hat{\mu} - 2g)\hat{\mu} + 2g^2} \leq \frac{(2\hat{\mu} + \mu)N}{\hat{\mu}} \sqrt{(\hat{\mu} - 2g)\hat{\mu} + 2g^2} \tag{26}$$

since $\|\hat{b}\| = N\sqrt{(\hat{\mu} - 2g)\hat{\mu} + 2g^2}$ and $\|\tilde{b}\| = \frac{|\mu - \hat{\mu}|N}{\hat{\mu}}\sqrt{(\hat{\mu} - 2g)\hat{\mu} + 2g^2}$. Note that

$\rho_0 := \rho - \sup_{t \in \mathbf{R}_{0+}} \|\tilde{\overline{A}}_0(t)\| > 0$ so that $(-\rho_0) < 0$ is larger than the maximum of the stability abscissas of

$\overline{A}(t)$ for $t \in \mathbf{R}_{0+}$ if Assertion 3 holds.

**Assertion 4**. If Assertion 3 holds then any solution of the forced system (17) to (22) satisfies the following inequality for some real constant $k_0 \geq 1$:

$$\|\overline{x}(t)\| \leq M(t) := k_0 e^{-\rho_0 t}\left(\|\overline{x}(0)\| + \|\overline{b}\|\int_0^t e^{\rho_0 \tau} d\tau\right) \to \frac{k_0}{\rho_0}\frac{(\hat{\mu} + |\mu - \hat{\mu}|)N}{\hat{\mu}}\sqrt{(\hat{\mu} - 2g)\hat{\mu} + 2g^2} \text{ as } t \to \infty$$

and the corresponding sub-states of $\overline{x}(t)$ satisfy:

$$\|\hat{x}(t)\| \leq \hat{M}(t) := k_0 e^{-\rho_0 t}\left(\|\hat{x}(0)\| + \|\hat{b}\|\int_0^t e^{\rho_0 \tau} d\tau\right) \to \frac{k_0}{\rho_0} N\sqrt{(\hat{\mu} - 2g)\hat{\mu} + 2g^2} \text{ as } t \to \infty$$

$$\|\tilde{x}(t)\| \leq \tilde{M}(t) := k_0 e^{-\rho_0 t}\left(\|\tilde{x}(0)\| + \|\tilde{b}\|\int_0^t e^{\rho_0 \tau} d\tau\right) \to \frac{k_0}{\rho_0}\frac{|\mu - \hat{\mu}|N}{\hat{\mu}}\sqrt{(\hat{\mu} - 2g)\hat{\mu} + 2g^2} \text{ as } t \to \infty$$

□

Note that $\|\tilde{x}(t)\| \to 0$ as $t \to \infty$ if $\hat{\mu} = \mu$ (and then $\|\hat{x}(t)\| \to \frac{k_0 N}{\rho_0}\sqrt{(\mu - 2g)\mu + 2g^2}$ as $t \to \infty$) and, if

in addition, $g = 0$ then $\|\hat{x}(t)\| \to \frac{k_0 \mu N}{\rho_0}$ as $t \to \infty$) or if $\hat{\mu} = g = 0$ (and then $\|\hat{x}(t)\| \to 0$ as $t \to \infty$).

Since $b := \left(\left(1 - \frac{g}{\hat{\mu}}\right)\mu N, 0, 0, \frac{g\mu N}{\hat{\mu}}\right)^T$ then $\|x(t)\| \leq \|M(t)\| \to \frac{k_0}{\rho_0}\frac{\hat{\mu} + |\mu - \hat{\mu}|N}{\hat{\mu}}\sqrt{(\hat{\mu} - 2g)\hat{\mu} + 2g^2}$ as

$t \to \infty$. However, this upper-bound can be improved if a version of Assertion 3 applied to the matrix function $A(t)$ leads to a smaller ratio of its corresponding constants than the ratio $k_0/\rho_0$ of the whole extended system.

**5. About the positivity of the SEIR extended epidemic model (17)-(22)**

Positive systems are those having nonnegative solutions in the sense that all the state components are nonnegative for all time [15-16]. Because of the nature of the SEIR epidemic model (1)-(4), it is required that it be a positive system for the implemented vaccination law. The extended SEIR- system has a unique solution for each initial state given by:

$$\overline{x}(t) = e^{\overline{A}_0 t}\left(\overline{x}(0) + \int_0^t e^{-\overline{A}_0 \tau}\left(\tilde{\overline{A}}_0(\tau)\overline{x}(\tau) + \overline{b}\right)d\tau\right) \tag{27}$$

From (10) and (13) the SEIR solution and its estimate through the observer are uniquely given by:

$$\hat{x}(t) = e^{\hat{A}_0 t}\left(\hat{x}(0) + \int_0^t e^{-\hat{A}_0 \tau}\left(\Delta\hat{A}(\tau)\hat{x}(\tau) + \hat{b}\right)d\tau\right) \tag{28}$$

$$x(t) = e^{A_0 t}\left(x(0) + \int_0^t e^{-A_0 \tau}\left[(\Delta A(\tau)x(\tau) + B(\tau)\hat{x}(\tau) + b)\right]d\tau\right) \tag{29}$$

In principle, it is apparently non necessary to require in addition that the estimation algorithm or the extended system be positive. However, note the following features by direct inspection of (28)-(29):

1) If $\hat{A}_0$ is a Metzler matrix, $\Delta\hat{A}(t) > 0; \forall t \in \mathbf{R}_{0+}$ and $\hat{b} > 0$ then $\hat{x}(0) \succ 0 \Leftrightarrow \hat{x}(t) \succ 0 ; \forall t \in \mathbf{R}_{0+}$.

2) If $\hat{x}(t) \succ 0 ; \forall t \in \mathbf{R}_{0+}$ (i.e. if $\hat{A}_0$ is a Metzler matrix, $\Delta\hat{A}(t) > 0; \forall t \in \mathbf{R}_{0+}$, $\hat{b} > 0$ and $\hat{x}(0) \succ 0$), $b > 0$, $\Delta A(t) > 0, B(t) > 0; \forall t \in \mathbf{R}_{0+}$ and $A_0$ is a Metzler matrix then $(x^T(0), \hat{x}^T(0))^T \succ 0 \Rightarrow x(t) \succ 0 ; \forall t \in \mathbf{R}_{0+}$.

3) $\hat{A}_0$ depends on the vaccination gain $k_1$, $\Delta\hat{A}(t)$ does not depend on $k_1$, but it depends on $k_5$ (see(21)), and $B(t)$ in (16) depends on $k_i$ ( i=1,…,5). Thus, if $k_1$ is chosen so that $\hat{A}_0$ is not a Metzler matrix (because its (4.1) entry is negative), $\hat{x}(t)$ may be non positive at some time for some nonnegative initial condition $\hat{x}(0)$. At the same time, the vector function $B(t)\hat{x}(t)$ may have some sufficiently negative component at some time "t" so that some corresponding component of $x(t)$ may be negative at that time. A close reasoning may be used for the case that $k_5$ is such that $B(t)$ is non positive (even if $\hat{A}_0$ is a Metzler matrix).

The following result follows from the above observations:

**Assertion 5**. The following properties hold:

(i) Assume that $A_0$ and $\hat{A}_0$ are Metzler matrices, $\Delta A(t) > 0, b + B(t)\hat{x}(t) \geq 0, \Delta\hat{A}(t) > 0; \forall t \in \mathbf{R}_{0+}$, $b > 0$ and $\hat{b} > 0$ then $(x^T(0), \hat{x}^T(0))^T \succ 0 \Rightarrow x(t) \succ 0, \hat{x}(t) \succ 0 ; \forall t \in \mathbf{R}_{0+}$. In other words, the extended system of state $(x^T(t), \hat{x}^T(t))^T$ is positive.

(ii) Assume that in Property (i) $\hat{A}_0$ fails to be a Metzler matrix because of the value $k_1$ in its (4,1) entry of $\Delta\hat{A}(t)$ is not positive due to the parameter $k_5$. Then, for initial conditions $\hat{x}(0) \succ 0$ which make $\hat{x}(t)$ to be non positive (i.e. with some negative component on a time interval), $x(t)$ can fail to be positive for all time for some $x(0) \succ 0$ and some such $k_1$ or $k_5$ with sufficiently large absolute values.

□

**Remarks 1**. It is required for modeling coherency that both the epidemic SEIR- model and its observer be positive dynamic systems. The condition of nonnegativeness of $(b + B(t)\hat{x}(t))$ for all time in Assertion 5 requires $g \geq 0$, and

$$(\hat{\mu} - g)N \geq (k_1 + k_5 \hat{I}(t))\hat{S}(t) + k_2 \hat{E}(t) + k_3 \hat{I}(t) + k_4 \hat{R}(t) \geq -gN ; \forall t \in \mathbf{R}_{0+}$$

which may be guaranteed by choosing the controller gains under the knowledge $N = \sum_{i=1}^{4} \hat{x}_i(t); \forall t \in \mathbf{R}_{0+}$

Also, $\hat{A}_0$ has to be a Metzler matrix, $\min_{t \in \mathbf{R}_{0+}} \Delta\hat{A}(t) > 0$ and $\hat{b} > 0$ ( Assertion 5 (i)) so that

$0 \leq g \leq \hat{\mu}$, $k_2 \leq 0, -\hat{\gamma} \leq k_3 \leq 0; 0 \leq k_4 \leq \hat{\omega}, k_1 + k_5 \hat{I}_r \geq 0$ in order that the observer be a positive system

$$\left(\hat{\beta}_1 + k_5\right)\left(\hat{I}_r - \max_{t \in \mathbf{R}_{0+}} \hat{I}(t)\right) \geq 0 \; ; \quad k_5\left(\min_{t \in \mathbf{R}_{0+}} \hat{I}(t) - \hat{I}_r\right) \geq 0 \tag{30}$$

This restricts the generality of the choice of the gains in the vaccination control law (9) since it would be needed to accomplish with the control gain constraint $0 \geq k_5 \geq -\beta_1$. However, if the requirement for the observer to be positive is removed then it is only needed that the SEIR- model (1)-(4) be positive under a modified vaccination law (9)

$$V(t) = \frac{1}{\hat{\mu} N}\left(k_1 \hat{S}(t) + k_2 \hat{E}(t) + k_3 \hat{I}(t) + k_4 \hat{R}(t) + k_5 \hat{S}(t)\hat{I}(t) + g N\right)$$

by requiring the weaker condition that $0 \leq g \leq \hat{\mu}$ and

$$\min(\sigma, \omega, \gamma) \geq 0 \; ; \quad I_r \geq \max_{t \in \mathbf{R}_{0+}} I(t) \tag{31} \quad \square$$

Note that while Assertion 5(i) is of sufficiency-type to guarantee positivity, the lack of all the joint above conditions in Assertion 5 (ii) refer to a necessary condition for positivity in such cases. Note also that positive and total population equal to N for all time implies necessary global stability so that we have directly the following:

**Assertion 6**. If Assertion 5 (i) holds then the extended SEIR –model (i.e. the combined SEIR- model plus its observer) is globally stable if all the initial populations and their estimates are nonnegative. Furthermore all the susceptible, infected, infectious and immune populations and their estimates are upper-bounded by N and the sum of all the populations and that of their estimates is equal to N at any time. The converse is not true, in general, so that if the extended SEIR – model is stable under Assertion 3 then such a model is not necessarily positive. $\square$

The positivity of the observer-based vaccination control is a nonnegative function if the conditions (30) hold, guaranteeing that the observer is a positive system, and furthermore the gain $0 \geq k_5 \geq \max(-\beta_1, -I_r / k_1) \geq -I(t)/k_1$. Thus, one has the following result:

**Assertion 7**. The vaccination control below is nonnegative for all time if $k_5 \geq 0$ and (30) hold:

$$V(t) = \frac{1}{\hat{\mu} N}\left(k_1 \hat{S}(t) + k_3 \hat{I}(t) + k_4 \hat{R}(t) + k_5 \hat{S}(t)\hat{I}(t) + g N\right) \tag{32}$$

With $-\hat{\gamma} \leq k_3 \leq 0$. Furthermore, the observer is a positive system under the vaccination control (32). If, in addition, (31) holds then the SEIR- model (1)-(4) is also a positive system under the vaccination control (32).

The more general vaccination control (9), namely

$$V(t) = \frac{1}{\hat{\mu} N}\left(k_1 \hat{S}(t) + k_2 \hat{E}(t) + k_3 \hat{I}(t) + k_4 \hat{R}(t) + k_5 \hat{S}(t)\hat{I}(t) + g N\right)$$

with $k_i \geq 0$; i=2,3,4, but not jointly zero, $k_5 \geq 0$ and subject to (30)-(31) is not guaranteed to be nonnegative for all time ( since the observer is not guaranteed to be a positive system). $\square$

The following vaccination nonnegative control combined of (9) and (32) may be used when the positivity of the observer is not imposed:

$$V(t) = \begin{cases} \overline{V}(t) & \text{if } 1 \geq \overline{V}(t) \geq 0 \\ \dfrac{1}{\hat{\mu} N}\left(k_1 \hat{S}(t) + k_4 \hat{R}(t) + k_5 \hat{S}(t)\hat{I}(t) + g N\right), & \text{otherwise} \end{cases} \qquad (33)$$

$$\overline{V}(t) := \dfrac{1}{\hat{\mu} N}\left(k_1 \hat{S}(t) + k_2 \hat{E}(t) + k_3 \hat{I}(t) + k_4 \hat{R}(t) + k_5 \hat{S}(t)\hat{I}(t) + g N\right) \qquad (34)$$

with $k_i \geq 0$; i =1,2,3,4,5 subject to (30)-(31). This nonnegative vaccination control keeps simultaneously as positive systems to both the SEIR-model and its observer.

### 6. Combined positivity, stability and tracking objective

The problem at hand requires as tracking objective simultaneously stability, positivity and tracking in the sense that the immune population converges to the total population what means that the total population including susceptible, infected and infectious populations is zero. Positivity implies global stability (while the converse is not always true) for all given combination of nonnegative initial populations but the suitable tracking objective is that the vaccination strategy lead asymptotically to a total immune population while the other populations converge to zero while keeping positivity and then global stability. This issue is now discussed through the design of the controller parameter g whose appropriate design is the basis for achievement of asymptotic tracking. Note that in order that the SEIR- model and its observer be both positive, it is needed that $0 \leq g \leq \hat{\mu}$ under Assertion 5 (i)- see Remark 1. The suitable tracking observer-based objective is $\hat{R}(t) \to N$ as $t \to \infty$ what implies $\hat{S}(t), \hat{E}(t), \hat{I}(t) \to 0$ as $t \to \infty$ with $\|\tilde{x}(t)\|$ being as small as possible as $t \to \infty$. Define

$$\hat{b}_1 := \hat{b}/N = (1,0,0,0)^T \hat{\mu} + (-1,0,0,1)^T g = e_1 \hat{\mu} + (e_4 - e_1) g \qquad (35)$$

where $e_i$ is the i-th unity vector in $\mathbf{R}^4$ with its i-th component being one. The above tracking objective may be achieved in some prescribed finite T or infinite time (i.e. asymptotically as $T \to \infty$). The asymptotic tracking objective is stated according to the relations of taking limits in (28) as $t \to \infty$ by designing a time-varying parameter $g \in PC^{(0)}(\mathbf{R}_{0+}; \mathbf{R})$ as follows:

$$g(t) = \begin{cases} \overline{g}(t) & \text{if } \overline{g}_M \geq \overline{g}(t) \geq 0 \text{ and } 0 \leq t \leq T \text{ or if } \hat{\mu} \geq \overline{g}(t) \geq 0 \text{ and } t > T \\ \overline{g}(f(t)), & \text{otherwise} \end{cases} \qquad (36)$$

$$\overline{g}(t) := \dfrac{1 - \int_0^t e_4^T e^{\hat{A}_0(t-\tau)}\left(e_1 \hat{\mu} + N^{-1} \Delta \hat{A}(\tau)\hat{x}(\tau)\right) d\tau}{\int_0^t e_4^T e^{\hat{A}_0(t-\tau)}(e_4 - e_1) d\tau} \qquad (37)$$

$$f(t) := max\left(t_1 < t : \overline{g}(t_1) \geq 0 \text{ and } g(t_1 + \tilde{t}) < 0; \forall \tilde{t} \in (0, t-t_1]\right) \quad \text{if } 0 \leq t \leq T \qquad (38)$$

$$f(t) := max\left(t_1 < t : \hat{\mu} \geq \overline{g}(t_1) \geq 0 \text{ and } g(t_1 + \tilde{t}) < 0; \forall \tilde{t} \in (0, t-t_1]\right) \quad \text{if } t > T \qquad (39)$$

with $\hat{\mu} \leq \bar{g}_M < \infty$ being a prefixed real constant which can be potentially large. It turns out the following main result:

**Theorem 1**. If $\hat{A}_0$ is a stability matrix then $g(t) \to g_\infty \in [0, \hat{\mu}]$. If the constant g is replaced with the piecewise continuous function defined by (36)-(39) in the SEIR- observer, (30) hold and the vaccination control (33)-(34) is used under non-negative initial conditions of the observer then the SEIR – observer is a positive system with all the populations being nonnegative and of total sum equal to N for all time while the vaccination control is also nonnegative for all time. Furthermore, $\hat{R}(t) \to N$; $\hat{S}(t), \hat{E}(t), \hat{I}(t) \to 0$ as $t \to \infty$.

If, in addition, (31) holds then the SEIR - model is also a positive system with populations nonnegative and equalizing N for all time. The estimation error fulfills the uniform upper-bound through time of Assertion 4.                                                                                          □

## 7. Simulation Examples

Consider the case when no vaccination is applied. The initial estimates of each population are $\hat{S}(0) = 250, \hat{E}(0) = 150, \hat{I}(0) = 150, \hat{R}(0) = 450$ individuals. Note that there is a great difference between the actual initial values and the estimated ones. Figure 1 shows the error between the actual and the estimated variables through time in the case when all the parameters of the model are known.

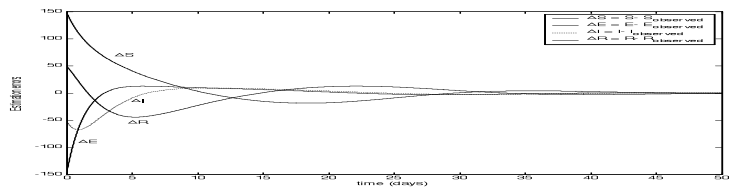

Figure 1. Observation error in the vaccination-free case and all the parameters known

Figure 2 shows, as a matter of example, the time-evolution of the susceptible and immune while Figure 3 shows the evolution of each population.

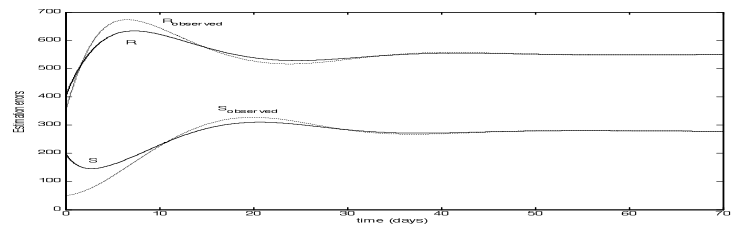

Figure 2. Evolution of the real and observed states for the susceptible and immune

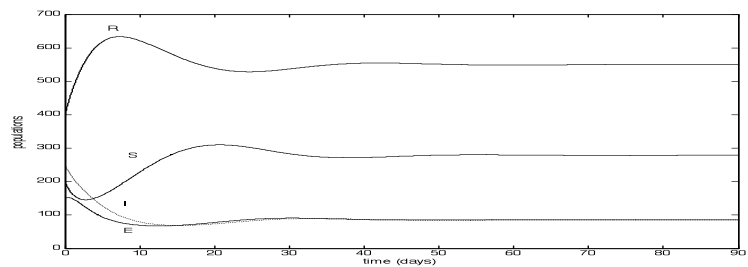

Figure 3. Evolution of the current populations through time

However, when the parameters of the system are not known and an estimation of its value is used, the observer does not track the real populations as shown in Figure 4 with estimated parameters

$$\frac{1}{\hat{\mu}} = 235 \text{ days}, \frac{1}{\hat{\sigma}} = 2 \text{ days}, \frac{1}{\hat{\omega}} = 14 \text{ days } \hat{\beta} = 1.46 days^{-1}, \hat{\gamma} = \hat{\sigma}.$$

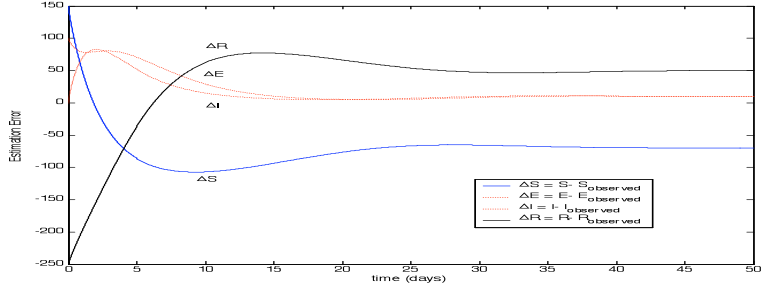

Figure 4. Observation error in the vaccination-free case and unknown model parameters

A vaccination strategy given by Equation (9) is now discussed. As commented in Section 5. In order to guarantee the positiveness of both systems, the control gains have to be chosen according Assertions 4 and 5 from Section 5 with $k_1 = 1, k_2 = -0.1 < 0, k_3 = -\hat{\gamma}, k_4 = 0.95\hat{\omega}, k_5 = -\hat{\beta}_1, g = \hat{\mu}$. Figures 5 and 6 show the time-evolution of the populations and the observation error, respectively.

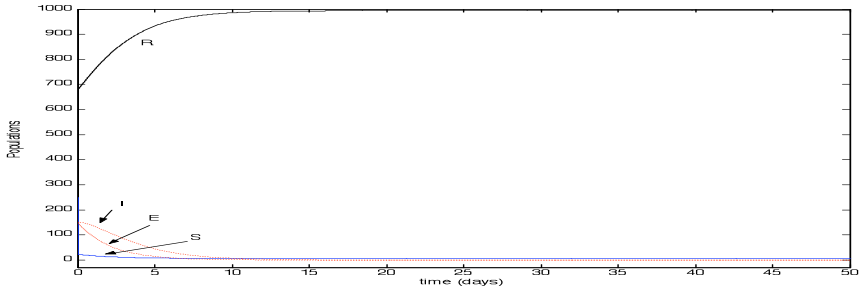

Figure 5. Evolution of the populations with vaccination and unknown parameters

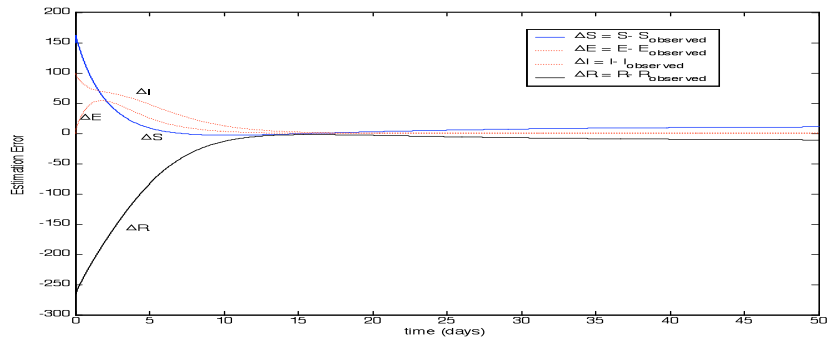

Figure 6. Observation error with vaccination and unknown parameters